# MODIFIED SMITH PREDICTOR FOR UNSTABLE LINEAR SYSTEMS

A. A. Pyrkin, K. Yu. Kalinin

ITMO University,

197101, St. Petersburg, Russia,

E-mail: pyrkin@itmo.ru

The paper presents a new control algorithm for unstable linear systems with input delay. In comparison with known analogues, the control law has been designed, which is a modification of the Smith predictor, and is the simplest one to implement without requiring complex integration methods. At the same time, the problem of stabilization of a closed system is effectively solved, ensuring the boundedness of all state variables and the exponential stability of the equilibrium point.

# МОДИФИЦИРОВАННЫЙ ПРЕДИКТОР СМИТА ДЛЯ НЕУСТОЙЧИВЫХ ЛИНЕЙНЫХ СИСТЕМ

А. А. Пыркин, К. Ю. Калинин

федеральное государственное автономное образовательное учреждение высшего образования «Национальный исследовательский университет ИТМО», 197101, г. Санкт-Петербург, Россия,
E-mail: pyrkin@itmo.ru

В работе представлен новый алгоритм управления для неустойчивых линейных систем с входным запаздыванием. В отличие от известных аналогов синтезирован закон управления, представляющий собой модификацию предиктора Смита, является наиболее простой в реализации, не требуя сложных методов интегрирования. При этом достаточно эффективно решена проблема стабилизации замкнутой системы, обеспечивая ограниченность всех переменных состояния и экспоненциальную устойчивость положения равновесия.

**Ключевые слова:** запаздывающее управление, предиктор Смита, управление по состоянию, неустойчивые системы.

## Введение

Одной из ключевых проблем как в классической, так и современной теории автоматического управления является синтез регуляторов для систем с запаздыванием. Причин возникновения запаздывания может быть несколько: удаленность объекта управления от системы управления, цифровые каналы передачи управляющих сигналов, конструктивные особенности объекта и другие.

Наличие временной задержки в контуре управления оказывает негативный эффект на свойства устойчивости замкнутой системы, и при не самых больших значениях запаздывания система теряет устойчивость. Физический смысл этого негативного влияния легко объясняется в терминах запаса устойчивости по фазе. Первая фундаментальная работа, в которой была изучена эта проблема и показано существование максимального запаздывания в контуре управления, соответствующего границе устойчивости, вышла почти сто лет назад [1]. Не менее важным прорывом стал Предиктор Смита [2, 3] – алгоритм управления, позволяющий при некоторых допущениях исключить влияние запаздывания на устойчивость замкнутой системы. Благодаря относительно простой реализации и убедительной эффективности этот подход широко распространился на практике. Тем не менее эти результаты позволяют синтезировать регуляторы только для устойчивых линейных систем с известными параметрами математической модели объекта управления, что существенно ограничивает область применения.

Следующим значительным прорывом стали работы [4-6], где были получены алгоритмы управления с предикцией переменных состояния для формирования стабилизирующей обратной связи для неустойчивых систем. Но по-прежнему рассматривались линейные системы с известными параметрами. Долгое время этот подход был и остается базовым для синтеза модифицированных версий предиктора, в том числе для систем с неизвестными параметрами, для нелинейных систем, для систем с распределенными параметрами, моделируемых уравнениями в частных производных [7, 8]. Ключевым недостатком этого подхода является использование неустойчивых динамических систем, включаемых в контур управления, и возникновение неустойчивой нуль-динамики в замкнутой системе. Было получено большое количество способов приближенного вычисления управляющего сигнала с целью сохранить устойчивость замкнутой системы, но даже такие подходы ограничены в своей применимости и могут рассматриваться как неконструктивные.

Крайне важным результатом по праву следует считать работу [9], в которой изучена проблема скрытой неустойчивой динамики, синтезирован наблюдатель переменных состояния этой нуль-динамики и получен робастный закон управления, позволяющий стабилизировать неустойчивые объекты с большим входным запаздыванием. В отличие от аналогов получена реализуемая схема синтеза регулятора, позволяющая снять ключевой недостаток предиктора [4]. С теоретической точки зрения важность этой работы сложно недооценить: решена почти полувековая нерешенная проблема неустойчивой скрытой динамики в предикторах для неустойчивых систем. Однако, стоит отметить, что структура закона управления не является очевидной, имеет в своей структуре достаточно много вспомогательных вычислений и переключений, что осложняет как реализацию, так и инженерное распространение этого решения.

В настоящей работе предложен новый алгоритм управления, отличающийся от [9] тем, что он не является модификацией [4-6], а представляющий собой самостоятельный и достаточно консервативный подход: добавление в Предиктор Смита корректирующего слагаемого, позволяющего стабилизировать замкнутую систему, в том числе и нуль-динамику. Как и в работе [9] в новом регуляторе используется прием сброса значений интегратора, однако, техническая реализация закона управления существенно проще: получен модифицированный предиктор Смита с корректирующим членом. Теоретическое доказательство устойчивости замкнутой системы также является новым и крайне перспективным для развития полученного решения и обобщения в будущем на нелинейные и адаптивные системы.

## Постановка задачи

Рассмотрим линейный объект с запаздыванием в канале управления
$$\dot{x}(t) = Ax(t) + Bu(t-D), \qquad (1)$$
где $x(t) \in R^n$ – измеряемый вектор переменных состояния, $A$, $B$ – полностью управляемая пара матриц с известными параметрами, $D \in R_+$ – известное постоянное запаздывание.

Требуется разработать закон управления $u(t)$, обеспечивающий асимптотическую устойчивость нулевого положения равновесия объекта $x = 0$ и ограниченность всех переменных состояния замкнутой системы.

## Синтез закона управления с предиктором

Выберем закон управления на основе классического Предиктора Смита
$$u(t) = Kx(t) + K\psi(t), \qquad (2)$$
$$\dot{\psi}(t) = A\psi(t) + Bu(t) - Bu(t-D) + L\zeta(t), \qquad (3)$$
с добавлением корректирующего члена $\zeta(t)$:
$$\zeta(t) = e^{AD}[x(t) - x(t-D) - \psi(t-D) + \varepsilon(t-D)] - \varepsilon(t), \qquad (4)$$
$$\dot{\varepsilon}(t) = A\varepsilon(t) + L\zeta(t), \quad \varepsilon(mT) = 0, \; m = 0,1,2,3,\ldots \qquad (5)$$
где матрицы $K$ и $L$ такие, что матрицы $F = A + BK$ и $H = A - L$ гурвицевы, а параметр $T$ будет определен позднее при анализе.

Заметим, что сигнал $\zeta(t)$ не является линейной комбинацией переменных состояния объекта $x(t)$ и регулятора $\psi(t)$, $\varepsilon(t)$, а зависит также от запаздывающих сигналов $x(t-D)$, $\psi(t-D)$ и $\varepsilon(t-D)$. Следовательно, $\zeta(t)$ может рассматриваться как дополнительная пространственная переменная, расширяющая динамику регулятора и замкнутой системы в целом.

Для производной $\zeta(t)$ справедлива модель:
$$\dot{\zeta}(t) = e^{AD}[\dot{x}(t) - \dot{x}(t-D) - \dot{\psi}(t-D) + \dot{\varepsilon}(t-D)] - \dot{\varepsilon}(t) =$$
$$= e^{AD}[Ax(t) + Bu(t-D) - Ax(t-D) - Bu(t-2D)] -$$
$$- e^{AD}[A\psi(t-D) + Bu(t-D) - Bu(t-2D) + L\zeta(t)] +$$

$$+e^{AD}[A\varepsilon(t-D)+L\zeta(t-D)]-A\varepsilon(t)-L\zeta(t)=$$
$$=e^{AD}[Ax(t)-Ax(t-D)-A\psi(t-D)+A\varepsilon(t-D)]-A\varepsilon(t)-L\zeta(t)=$$
$$=Ae^{AD}[x(t)-x(t-D)-\psi(t-D)+\varepsilon(t-D)]-A\varepsilon(t)-L\zeta(t)=$$
$$=A\zeta(t)-L\zeta(t)=$$
$$=(A-L)\zeta(t). \tag{6}$$

Отметим, что модель (6) глобально экспоненциально устойчива в силу гурвицевости матрицы $A-L$, однако, в дискретные моменты времени $mT$ и $mT+D$ значение переменной $\zeta(t)$ скачкообразно меняется согласно уравнениям (4) и (5).

Введем в рассмотрение замену координат:
$$z(t)=\zeta(t)+\varepsilon(t), \tag{7}$$

Вычислим производную
$$\dot{z}(t)=(A-L)\zeta(t)+A\varepsilon(t)+L\zeta(t)=$$
$$=Az(t),$$

На основе (4) получим выражение для запаздывающего управления
$$u(t-D)=Kx(t-D)+K\psi(t-D)=$$
$$=Kx(t)+K\varepsilon(t-D)-Ke^{-AD}z(t)$$

и перепишем модель для переменных $x(t)$, $\psi(t)$ и $\varepsilon(t)$:
$$\dot{x}(t)=Ax(t)+BKx(t)-BKe^{-AD}z(t)+BK\varepsilon(t-D)=$$
$$=Fx(t)+BK\varepsilon(t-D)-BKe^{-AD}z(t),$$
$$\dot{\psi}(t)=A\psi(t)+BKx(t)+BK\psi(t)-$$
$$-BKx(t)-BK\varepsilon(t-D)+BKe^{-AD}z(t)+L\zeta(t)=$$
$$=F\psi(t)+BKe^{-AD}z(t)+Lz(t)-BK\varepsilon(t-D)-L\varepsilon(t),$$
$$\dot{\varepsilon}(t)=A\varepsilon(t)+Lz(t)-L\varepsilon(t)=$$
$$=H\varepsilon(t)+Lz(t).$$

Перепишем модель замкнутой системы (9)-(13) в компактном виде:
$$\dot{x}(t)=Fx(t)+BK\varepsilon(t-D)-BKe^{-AD}z(t), \tag{8}$$
$$\dot{\psi}(t)=F\psi(t)+BKe^{-AD}z(t)+Lz(t)-BK\varepsilon(t-D)-L\varepsilon(t), \tag{9}$$

$$\dot{\varepsilon}(t) = H\varepsilon(t) + Lz(t). \tag{10}$$

$$\dot{z}(t) = Az(t). \tag{11}$$

Нетрудно видеть, что динамика переменной $z(t)$ может быть неустойчива в силу свойств матрицы $A$, которая по постановке задачи не обязательно гурвицева. Если не использовать корректировку переменной $\varepsilon(t)$, то замкнутая система может быть неустойчивой, поскольку на вход устойчивых по входу подсистем (8)-(10) попадает неограниченный сигнал $z(t)$. Далее будем рассматривать замкнутую систему с принудительным обнулением переменной $\varepsilon(t)$ с периодом $T$.

## Основной результат

Заметим, что функция времени

$$z(t) = e^{AD}[x(t) - x(t - D) - \psi(t - D) + \varepsilon(t - D)] \tag{12}$$

не является непрерывной поскольку она алгебраически зависит от функции $\varepsilon(t - D)$, имеющей разрывы первого рода. Причем значение переменной $z(t)$ будет меняться скачком в моменты времени $mT + D$. Для переменных $x(t)$ и $\psi(t)$ можно показать их непрерывность. Далее необходимо проанализировать, при каких условиях сброса переменной $\varepsilon(t)$ замкнутая система является асимптотически устойчивой.

В силу (11) между переключениями в моменты времени $[t_k, t_{k+1})$ для переменной $z(t)$ справедливо выражение

$$z(t) = e^{A(t - t_k)} z(t_k). \tag{13}$$

Исследуем последовательность значений функции $z(t)$ для моментов времени $t = mT + D$ для $m \in Z_{\geq 0}$

$$z(mT + D) = e^{AD}[x(mT + D) - x(mT) - \psi(mT)], \tag{14}$$

где уже учтено тождество $\xi(t - D) = 0, \forall t = mT + D$.

Для этого рассмотрим вспомогательную непрерывную функцию времени

$$\xi(t) = x(t) - x(t - D) - \psi(t - D) = e^{-AD}z(t) - \varepsilon(t - D) \tag{15}$$

и ее модель

$$\dot{\xi}(t) = e^{-AD}\dot{z}(t) - \dot{\varepsilon}(t-D) =$$
$$= Ae^{-AD}z(t) - A\varepsilon(t-D) - L\zeta(t-D) =$$
$$= A\xi(t) - Lz(t-D) + L\varepsilon(t-D) =$$
$$= A\xi(t) - Lz(t-D) + Le^{-AD}z(t) - L\xi(t) =$$
$$= H\xi(t) + Le^{-AD}z(t) - Lz(t-D). \qquad (16)$$

где заметим, что $\xi(0) = x(0)$ и $\xi(D) = x(D) - x(0) - \psi(0)$.

Интегрируя (16), получим выражение

$$\xi(t) = e^{Ht}x(0) + \int_0^t e^{H(t-s)}L(e^{-AD}z(s) - z(s-D))ds. \qquad (17)$$

Далее запишем выражение для последовательности $\xi_m = \xi(mT + D)$:

$$\xi_m = e^{H(mT+D)}x(0) + \int_0^{mT+D} e^{H(mT+D-s)}L(e^{-AD}z(s) - z(s-D))ds. \qquad (18)$$

Вычислим значение последовательности на следующем шаге

$$\xi_{m+1} = e^{H(mT+T+D)}x(0) + \int_0^{mT+T+D} e^{H(mT+T+D-s)}L(e^{-AD}z(s) - z(s-D))ds.$$

Подставим в последнее уравнение выражение $e^{H(mT+D)}x(0)$ из (18):

$$\xi_{m+1} = e^{HT}\xi_m - e^{HT}\int_0^{mT+D} e^{H(mT+D-s)}L\big(e^{-AD}z(s) - z(s-D)\big)ds +$$
$$+ e^{HT}\int_0^{mT+T+D} e^{H(mT+D-s)}L\big(e^{-AD}z(s) - z(s-D)\big)ds =$$
$$= e^{HT}\xi_m + e^{HT}\int_{mT+D}^{mT+T+D} e^{H(mT+D-s)}L\big(e^{-AD}z(s) - z(s-D)\big)ds =$$
$$= e^{HT}\xi_m + e^{HT}\int_{mT+D}^{mT+2D} e^{H(mT+D-s)}L\big(e^{-AD}z(s) - z(s-D)\big)ds +$$
$$+ e^{HT}\int_{mT+2D}^{mT+T+D} e^{H(mT+D-s)}L\big(e^{-AD}z(s) - z(s-D)\big)ds.$$

В моменты времени $t = mT + 2D$ скачком меняется значение функции $z(t - D)$, поэтому для вычисления интеграла необходимо разбить его на два интервала $[mT + D; mT + 2D)$ и $[mT + 2D; mT + D + T)$ и считать отдельно.

Для первого интервала значение функции $z(t - D)$ равно

$$z(t-D) = e^{A(t-(m-1)T-2D)}z((m-1)T + D), \quad mT + D \leq t < mT + 2D,$$

что соответствует непрерывному росту с момента предыдущей коррекции $\varepsilon(t)$ на шаге $m - 1$. В момент времени $mT + 2D$ функция $z(t - D)$ меняет скачком свое значение, что соответствует шагу $m$. Затем на втором интервале имеем

$$z(t - D) = e^{A(t-mT-2D)}z(mT + D), \quad mT + 2D \le t < mT + D + T.$$

В течение обоих интервалов в период $[mT + D; mT + D + T)$ функция z(t) имеет вид

$$z(t) = e^{A(t-mT-D)}z(mT + D).$$

Далее заметим, что

$$e^{-AD}z(t) - z(t - D) \equiv 0, \quad mT + 2D \le t < mT + D + T,$$

откуда видим, что интеграл на втором интервале равен 0. На первом интервале $mT + D \le t < mT + 2D$ вычислим:

$$e^{-AD}z(t) - z(t - D) =$$
$$= e^{-AD}e^{A(t-mT-D)}z(mT + D) - e^{A(t-(m-1)T-2D)}z((m-1)T + D) =$$
$$= e^{A(t-mT-2D+T)}[e^{-AT}z(mT + D) - z(mT + D - T)]$$

Тогда

$$\xi_{m+1} = e^{HT}\xi_m + e^{HT}\int_{mT+D}^{mT+2D} e^{H(mT+D-s)}L\big(e^{-AD}z(s) - z(s - D)\big)ds =$$
$$= e^{HT}\xi_m +$$
$$+ e^{HT}\int_{mT+D}^{mT+2D} e^{H(mT+D-s)}Le^{A(s-mT-2D+T)}[e^{-AT}z(mT + D) - z(mT + D - T)]ds$$

Заметим, что в силу (15) и правила обнуления (5) переменной $\varepsilon(t)$, справедливы соотношения

$$z(mT + D) = e^{AD}\xi_m, \quad z(mT + D - T) = e^{AD}\xi_{m-1},$$

тогда получим выражение

$$\xi_{m+1} = e^{HT}\xi_m + e^{HT}\int_{mT+D}^{mT+2D} e^{H(mT+D-s)}Le^{A(s-mT-2D+T)}e^{AD}[e^{-AT}\xi_m - \xi_{m-1}]ds =$$
$$= e^{HT}\xi_m + e^{H(mT+D+T)}\int_{mT+D}^{mT+2D} e^{-Hs}Le^{As}ds \times e^{A(-mT-D+T)}[e^{-AT}\xi_m - \xi_{m-1}].$$

**Лемма 1.** Справедливо соотношение для интеграла
$$\int e^{-Hs} L e^{As} ds = -e^{-Hs} e^{As} + \text{const.}$$

**Доказательство леммы.**

Дифференцируя функцию $e^{-Hs} e^{As}$, нетрудно видеть
$$\frac{d}{dt}(-e^{-Hs} e^{As}) = e^{-Hs} H e^{As} - e^{-Hs} e^{As} A =$$
$$= e^{-Hs}(A - L) e^{As} - e^{-Hs} A e^{As} =$$
$$= e^{-Hs} L e^{As},$$

что соответствует подынтегральному выражению, что и требовалось доказать. ∎

Тогда справедливо соотношение:
$$\xi_{m+1} = e^{HT} \xi_m + W(T)[e^{-AT} \xi_m - \xi_{m-1}],$$
где
$$W(T) = e^{H(mT+D+T)}\bigl(-e^{-H(mT+2D)} e^{A(mT+2D)} + e^{-H(mT+D)} e^{A(mT+D)}\bigr) e^{A(-mT-D+T)}$$
$$= -e^{HT-HD} e^{AD+AT} + e^{HT} e^{AT},$$

что позволяет записать характеристическое уравнение для последовательности $\xi_m$:
$$\xi_{m+1} - (e^{HT-HD} e^{AD}) \xi_m + (e^{HT-HD} e^{AD+AT} - e^{HT} e^{AT}) \xi_{m-1} = 0 \qquad (19)$$
или в блочном матричном виде:
$$\begin{pmatrix} \xi_m \\ \xi_{m+1} \end{pmatrix} = \begin{pmatrix} 0 & I \\ -G_2 & -G_1 \end{pmatrix} \begin{pmatrix} \xi_{m-1} \\ \xi_m \end{pmatrix}, \qquad (20)$$
с единичной матрицей $I$ размерности $n$ и матричными коэффициентами
$$G_1 = -e^{HT-HD} e^{AD}, \quad G_2 = e^{HT}(e^{-HD} e^{AD} - I) e^{AT}.$$

**Лемма 2.** Существует $T_0 > D$ такое, что $\forall T \geq T_0$ последовательность $\xi_m$ экспоненциально сходится к $0$.

**Доказательство леммы.** Рассмотрим матрицу $P = \begin{pmatrix} I & 0 \\ 0 & 2I \end{pmatrix}$ и функцию Ляпунова $V(m) = \chi^\top(m) P \chi(m)$ с вектором состояния $\chi(m) = \mathrm{col}(\xi_{m-1}, \xi_m)$. Тогда

$$V(m) - V(m+1) = \chi^\top(m) P \chi(m) - \chi^\top(m+1) P \chi(m+1) =$$

$$= \chi^\top(m) \left[ \begin{pmatrix} I & 0 \\ 0 & 2I \end{pmatrix} - \begin{pmatrix} 0 & -G_2^\top \\ I & -G_1^\top \end{pmatrix} \begin{pmatrix} I & 0 \\ 0 & 2I \end{pmatrix} \begin{pmatrix} 0 & I \\ -G_2 & -G_1 \end{pmatrix} \right] \chi(m) =$$

$$= \chi^\top(m) \left[ \begin{pmatrix} I & 0 \\ 0 & I \end{pmatrix} - 2 \begin{pmatrix} G_2^\top G_2 & G_2^\top G_1 \\ G_1^\top G_2 & G_1^\top G_1 \end{pmatrix} \right] \chi(m).$$

Поскольку матрица $H = A - L$ гурвицева в силу соответствующего выбора $L$, то

$$\lim_{T \to \infty} e^{HT} = 0.$$

Более того, коэффициенты $L$ могут быть выбраны так, чтобы выполнялось

$$\mathbf{Re}\{\lambda_{\max}\{H\} + \lambda_{\max}\{A\}\} < 0,$$

что будет гарантировать затухание функции $e^{HT}$ быстрее возможного роста функции $e^{AT}$ для произвольной матрицы $A$, обеспечивая

$$\lim_{T \to \infty} e^{HT} e^{AT} = 0.$$

Поскольку функция $N(T) = 2 \begin{pmatrix} G_2^\top(T) G_2(T) & G_2^\top(T) G_1(T) \\ G_1^\top(T) G_2(T) & G_1^\top(T) G_1(T) \end{pmatrix}$ является непрерывной, ограниченной для $\forall T \geq 0$, а все ее элементы с ростом $T$ стремятся к нулю $\lim_{T \to \infty} N(T) = 0$, то существует $T_0 > D$ такое, что $\forall T \geq T_0$ справедливо неравенство

$$\|N(T)\| \leq \alpha I, \quad 0 < \alpha < 1.$$

Следовательно,
$$V(m) - V(m+1) = \chi^\top(m)[I_{2n} - N(T)]\chi(m) > (1-\alpha)\chi^\top(m)\chi(m) \geq \frac{1-\alpha}{2} V(m).$$

Таким образом, при условии $T \geq T_0$

$$V(m+1) < \beta V(m) < \beta^m V(0),$$

откуда следует равномерная экспоненциальная сходимость последовательности $V(m)$ к нулю с коэффициентом $\beta = \frac{1+\alpha}{2} < 1$ и экспоненциальная сходимость последовательности $\xi_m$, что и требовалось доказать. ∎

Напомним, что последовательность $\xi_m$ представляет собой значения непрерывной функции $\xi(t)$ в моменты времени $t = mT + D$, и она сходится к нулю. Также заметим, что в силу определения последовательность $\xi_m$ соответствует значениям функции $z(t)$ в моменты времени $t = mT + D$. Функция $z(t)$ на интервале $[mT + D; mT + D + T)$ изменяется согласно выражению (13). Максимальное значение $z(t)$ на интервале между переключениями $[mT + D; mT + D + T)$ можно оценить:

$$\|z(t)\| \leq \max_{0 \leq s < T} \|e^{As}\| \cdot \|z(mT + D)\|.$$

Так как последовательность $z(mT + D)$ экспоненциально стремится к нулю с увеличением $m$, то можно показать экспоненциальную сходимость к нулю всех элементов вектора $z(t)$, указав соответствующую мажоранту.

**Утверждение.** Объект управления с запаздывающим управлением (1) и регулятором на основе предиктора Смита (2), (3) с корректирующим членом (4), (5) обеспечивает глобальную экспоненциальную устойчивость нулевого положения равновесия замкнутой системы.

**Доказательство утверждения.** Замкнутая система описывается моделью (8)-(11) с переключениями сигналов $\varepsilon(t)$ по правилу (5) и сигнала $z(t)$ в соответствие с (14). Все элементы вектора $z(t)$ экспоненциально сходятся к нулю, поскольку норма этого вектора ограничена мажорирующей затухающей экспонентой. Далее нетрудно видеть, что сигнал $\varepsilon(t)$ также сходится экспоненциально к нулю в силу модели (10). Аналогично можно показать ограниченность и экспоненциальную сходимость к нулю всех переменных состояния модели (8)-(11), откуда следует экспоненциальная устойчивость нулевого положения равновесия замкнутой системы. ∎

**Пример численного моделирования**

Рассмотрим неустойчивый объект управления (1) с параметрами $A = \begin{bmatrix} 0 & 1 \\ 0{,}1 & 0 \end{bmatrix}$, $B = \begin{bmatrix} 0 \\ 1 \end{bmatrix}$, $x(0) = \begin{bmatrix} -1 \\ 1 \end{bmatrix}$ и регулятор (2)-(5) с параметрами $K = \begin{bmatrix} -20 & -30 \end{bmatrix}$, $L = \begin{bmatrix} 2 & 0{,}5 \\ 3 & 0 \end{bmatrix}$. На рисунках 1, 2 представлены результаты моделирования для различных значений запаздывания $D$ и интервала сброса интегратора $T$ с использованием метода интегрирования **ode1 (Euler)** с фиксированным шагом $10^{-4}$ с.

**Заключение**

В работе представлен принципиально новый предиктор для неустойчивых линейных систем с входным запаздыванием, отличающийся более простой структурой в реализации, строгим аналитическим доказательством устойчивости. Разработанный подход открывает широкие возможности для дальнейших обобщающих работ, в которых можно синтезировать предиктор по выходу, для нелинейных систем, для систем с неизвестными параметрами. Консервативная структура как регулятора, так и доказательства его эффективности позволяет утверждать о том, что получена база для нового фундаментального метода в теории автоматического управления.

СПИСОК ЛИТЕРАТУРЫ

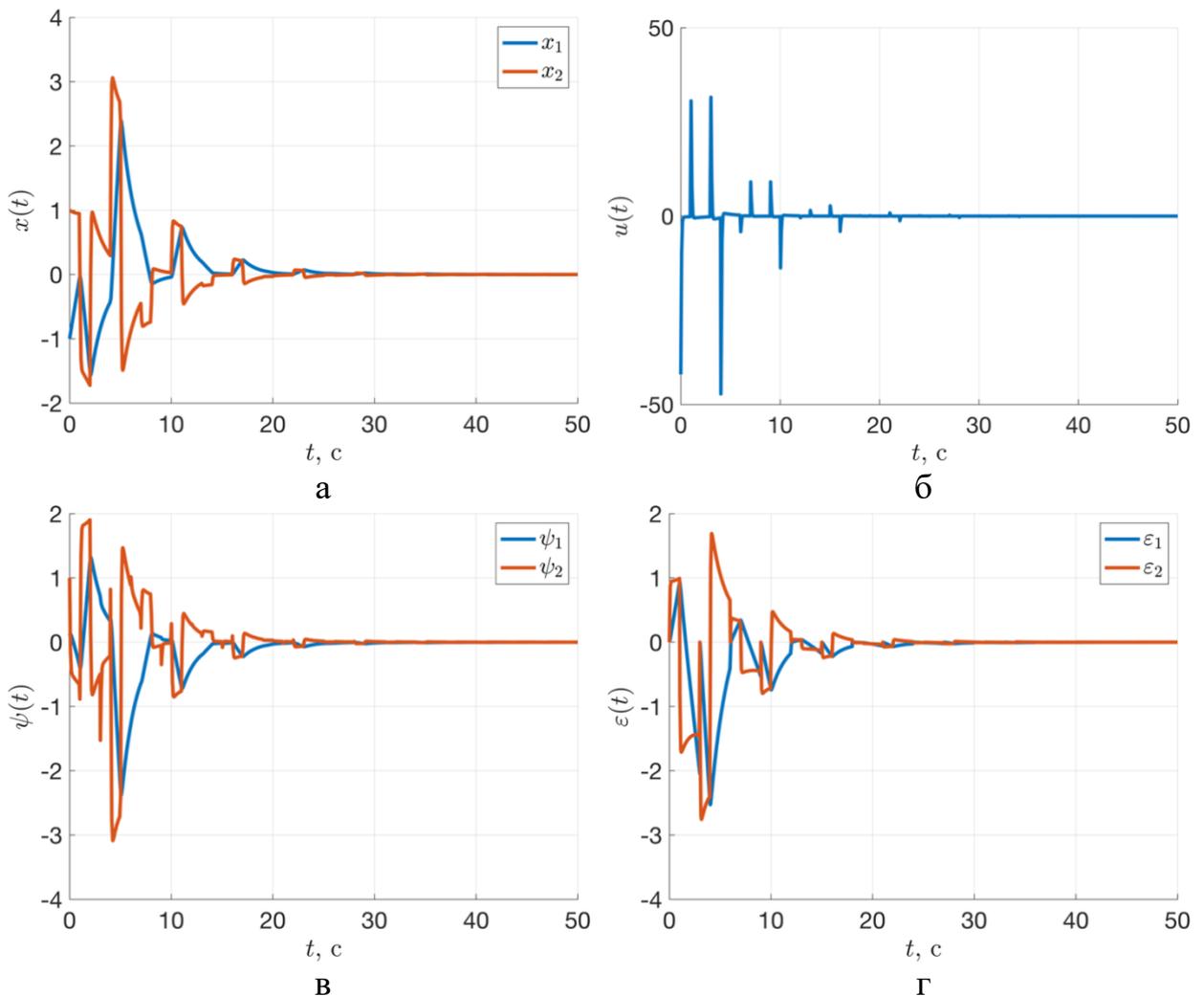

Рисунок 1 – Переходные процессы для $D = 1$ и $T = 5$

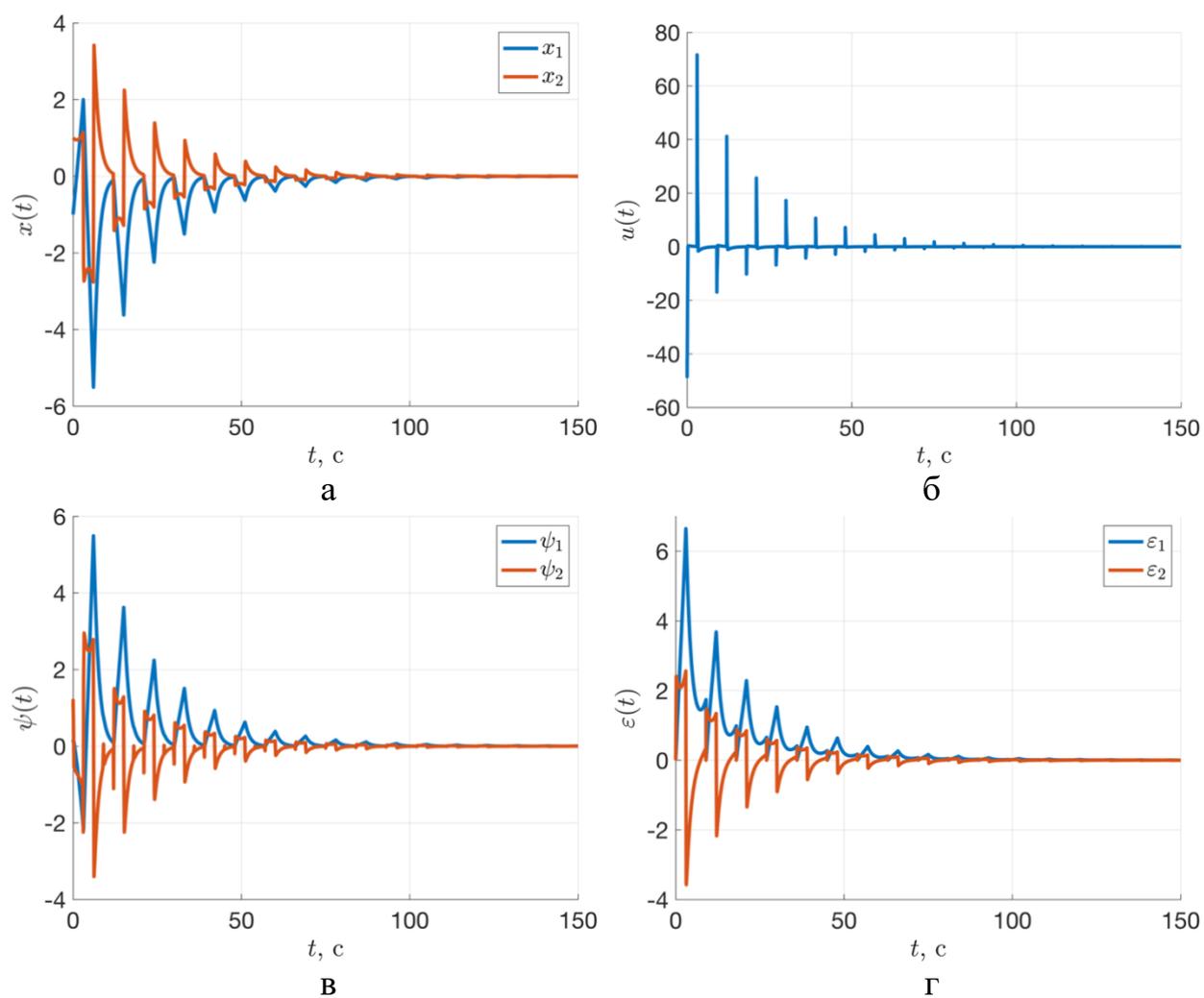

Рисунок 2 – Переходные процессы для $D = 3$ и $T = 30$